\documentclass[12pt]{article}

\usepackage{latexsym}
\usepackage{amsmath}
\usepackage{amsfonts}
\usepackage{amssymb}

\newcommand{\bean}{\begin{eqnarray}}
\newcommand{\eean}{\end{eqnarray}}
\newcommand{\bea}{\begin{eqnarray*}}
\newcommand{\eea}{\end{eqnarray*}}
\newcommand{\bsa}{\begin{subarray}{c}}
\newcommand{\esa}{\end{subarray}}
\newcommand{\bi}{\begin{itemize}}
\newcommand{\ei}{\end{itemize}}

\newtheorem{lemma}{Lemma}[section]

\newtheorem{propn}[lemma]{Proposition}

\begin{document}

\title{\bf On the Fourier coefficients of 2-dimensional vector-valued modular forms}
\author{Geoffrey Mason\thanks{Supported by NSF and NSA.} \\
University of California, Santa Cruz}
\date{}
\maketitle
 \begin{abstract}
\noindent
 Let $\rho: SL(2, \mathbb{Z}) \rightarrow GL(2, \mathbb{C})$ be an irreducible representation of the modular group such that $\rho(T)$ has finite order $N$. We study holomorphic vector-valued modular forms $F(\tau)$ of integral weight associated to $\rho$ which have \emph{rational} Fourier coefficients. (These span the complex space of all integral weight vector-valued modular forms associated to $\rho$.) As a special case of the main Theorem, we prove that if
 $N$ does \emph{not} divide $120$ then every nonzero $F(\tau)$ has Fourier coefficients with
 \emph{unbounded denominators}.
\end{abstract}

\section{Introduction}
 Let $\Gamma = SL(2, \mathbb{Z})$ be the inhomogeneous modular group with standard generators
 \begin{eqnarray*}
S = \left(\begin{array}{cc}0 & -1 \\1 & 0\end{array}\right), 
\ T = \left(\begin{array}{cc}1 & 1 \\0 & 1\end{array}\right).
\end{eqnarray*}
We will be dealing with $2$-dimensional irreducible representations
$\rho: \Gamma \rightarrow GL(2, \mathbb{C})$
 such that $\rho(T)$ is a diagonal \emph{unitary} matrix
 \begin{eqnarray}\label{rhoTdef}
\rho(T) = \left(\begin{array}{cc} e^{2\pi i m_1}&0 \\ 0&e^{2\pi i m_2}\end{array}\right), \ \ \ (1 > m_1
> m_2 \geq 0).
\end{eqnarray}

  A holomorphic vector-valued modular form 
 of integral weight $k$ associated to $\rho$ consists of the following data: \\
 (i) a pair of holomorphic functions 
 $f_i(\tau), \ i=1, 2,$ defined on the complex upper half-plane $\frak{H}$
  with  $q$-expansions
 \begin{eqnarray}\label{qexp0}
f_i(\tau) = \sum_{n \geq 0} a_{n, i}q^{m_i+n},
\end{eqnarray}
(ii) functional equations
\begin{eqnarray}\label{fe0}
F|_k \gamma (\tau) = \rho(\gamma)F(\tau), \ \ (\gamma \in \Gamma).
\end{eqnarray}
Here, $F(\tau) = (f_1(\tau), f_2(\tau))^t$ is\footnote{superscript $t$ means \emph{transpose}} the column vector whose components are the functions
$f_i(\tau)$ and $|_k$ is the usual stroke operator in weight $k$ applied to each component. The coefficients $a_{n, i}$ in (\ref{qexp0}) are  the \emph{Fourier coefficients} of $f_i(\tau)$,
or $F(\tau).$  The set $\mathcal{H}(k, \rho)$ of all holomorphic vector-valued modular forms of weight $k$ associated to  $\rho$ is a finite-dimensional $\mathbb{C}$-linear space (\cite{KM}).

\medskip
We say that \emph{$\rho$ is modular} of \emph{level $N$} if ker$\rho$ is a congruence subgroup of level $N$, i.e. ker$\rho \supseteq \Gamma(N)$.  In this case, the component functions of a vector-valued modular form
$F(\tau) \in \mathcal{H}(k, \rho)$ are classical holomorphic modular forms of weight $k$ and level $N$. 
For the purposes of the present paper,
the \emph{projective level} of $\rho$ is a more useful invariant. If
$\overline{\rho}$ is the projective representation of $\Gamma$ defined by  the composition 
\begin{eqnarray}\label{projrep}
\Gamma \stackrel{\rho}{\rightarrow} GL(2, \mathbb{C}) \stackrel{\pi}{\rightarrow} PGL(2, \mathbb{C})
\end{eqnarray}
($\pi$ is canonical projection), we define the projective level of $\rho$ to
be the \emph{order} of 
$\overline{\rho}(T) \in PGL(2, \mathbb{C})$. We emphasize that we are \emph{not} assuming that ker$\bar{\rho}$ is a congruence subgroup.

\medskip

There are infinitely many equivalence classes of $2$-dimensional irreducible $\rho$, but it turns out 
(\cite{M2}) that just $54$ of them are modular. The ordinary level $N$ in these cases is some divisor of $120$, and  for such an $N$ some
classes of $\rho$ are modular, while others may not be. On the other hand, it can be deduced from the tables in \cite{M2} (cf. Theorem \ref{thm2.4} below) that
\begin{eqnarray}\label{ab1}
\rho \ \mbox{is modular} \Leftrightarrow M \leq 5. 
\end{eqnarray} (\ref{ab1}) is an \emph{algebraic} characterization of those $\rho$ which are modular.
The present paper is concerned with the problem of characterizing the modular $\rho$ by means of
\emph{arithmetic} properties of the Fourier coefficients of associated vector-valued modular forms.
We are mainly interested in the space $\mathcal{H}(k, \rho)_{\mathbb{Q}} \subseteq \mathcal{H}(k, \rho)$
consisting of the $F(\tau)$ whose Fourier coefficients lie in $\mathbb{Q}$. As we will see (Lemma \ref{lemma5.0}), $\mathcal{H}(k, \rho)_{\mathbb{Q}}$ is a $\mathbb{Q}$-form for
$\mathcal{H}(k, \rho)$ (i.e. $\mathcal{H}(k, \rho)= \mathbb{C} \otimes_{\mathbb{Q}}\mathcal{H}(k, \rho)_{\mathbb{Q}}$) if the projective level of $\rho$ is finite. 

\medskip
Suppose that $F(\tau) \in \mathcal{H}(k, \rho)_{\mathbb{Q}}$ has component functions (\ref{qexp0}).
We say that $F(\tau)$ has \emph{bounded denominators} if there is a nonzero integer $B$ such that $Ba_{n, i} \in \mathbb{Z} \ (n\geq 0,
i=1, 2)$. Otherwise, $F(\tau)$ has \emph{unbounded denominators}. We can now state the

\medskip
\noindent
Conjecture: Suppose that $\rho$ has finite projective level. Then $\rho$ is
modular if, and only if, there is some nonzero $F(\tau) \in \mathcal{H}(k, \rho)_{\mathbb{Q}}$
which has bounded denominators.

\medskip
We make several remarks. The irreducibility of $\rho$ is implicitly assumed in the Conjecture. If $\rho$ is modular then, as we have explained, the components of a vector-valued modular form
in $\mathcal{H}(k, \rho)$ are ordinary modular forms, in which case the bounded denominator property of Fourier coefficients is well-known. So the Conjecture really concerns the implication 
\emph{bounded denominators $\Rightarrow \rho$ modular}. 
It is one of a hierarchy of similar conjectures about the modularity of vector-valued modular forms (of arbitrary finite dimension) whose Fourier coefficents are algebraic with bounded denominators. Other special cases that have been discussed in the literature include modular forms on noncongruence subgroups (cf Atkin-Swinnerton-Dyer \cite{AS} and
Kurth-Long \cite{KL1}, \cite{KL2}),  and generalized modular forms (Kohnen-Mason (\cite{KoM1}, \cite{KoM2}).

\medskip
The main result of the present paper is the following:

\medskip
\noindent
{\bf Theorem 1}\label{thmmain1} Suppose that $\rho$ has finite projective level $M$, and that $M$ does \emph{not} divide $60$.
 Then the components of \emph{every} nonzero vector-valued modular form
$F(\tau) \in \mathcal{H}(k, \rho)_{\mathbb{Q}}$ have unbounded denominators.

\medskip
 By Theorem \ref{thmmain1}, a counterexample  to the Conjecture
  necessarily has projective level $M|60$. There are 
approximately $350$ equivalence classes of $\rho$ satisfying this  condition, including of course the $54$ classes which are modular. In particular, Theorem \ref{thmmain1} proves the Conjecture for all but finitely many equivalence classes of $\rho$. 
  
\bigskip
The proof of Theorem \ref{thmmain1} depends on results in \cite{M1} and \cite{M2} (see also \cite{MM})
describing $\mathcal{H}(\rho) = \oplus_k \mathcal{H}(k, \rho)$
as a module over a certain ring $\mathcal{R}$ of differential operators.
There is a minimal weight $k_0$ for which
$\mathcal{H}(k_0, \rho)$ is nonzero, and this space is $1$-dimensional with basis $F_0(\tau)$, say.
$F_0(\tau)$ \emph{generates} $\mathcal{H}(\rho)$ considered as $\mathcal{R}$-module, and the general idea of the proof is to reduce questions about arbitrary 
$F(\tau)$ to questions about $F_0(\tau)$. The components of $F_0(\tau)$ span the solution space of
a certain \emph{modular linear differential equation} (\cite{M1}) which has $q=0$ as a regular-singular point, and the  Fuchsian theory provides a recursive formula for the corresponding Fourier coefficients.
 Assuming (as we may) that $F_0(\tau)$ has rational Fourier coefficients and
 that $M$ does not divide $60$, we can exploit the recursion to obtain (Proposition \ref{thm1.4}) the \emph{exact} power of $p$ dividing the denominators of the Fourier coefficients 
 $a_{n, i}$ of $F_0(\tau)$ whenever $p$ is a prime dividing $M/(M, 60)$. The power of $p$  is \emph{strictly increasing} for 
 $n \rightarrow \infty$, and in particular
 $F_0(\tau)$ satisfies the unbounded denominator property.
Together with the structure of $\mathcal{H}(\rho)$ as $\mathcal{R}$-module, this result can then 
be used to deduce Theorem 1 for general $F(\tau)$.

\medskip
The paper is organized as follows. In Section 2 we cover the background needed from \cite{M1}, \cite{M2}.  In Section 3 we discuss the projective level
$M$ and related invariants, and in particular we give (Proposition \ref{thm2.4}) a direct proof of (\ref{ab1}) which does not rely on the tables in \cite{M2}. Section 4 contains the proof of Theorem 1, and Section 5 
contains some concluding remarks.
 
 \section {Background}
  We review notation and results we will need from \cite{M1}, \cite{M2}. The $\mathbb{Z}$-graded ring of holomorphic modular forms on $\Gamma$ is 
\begin{eqnarray*}
\mathcal{M} = \oplus_{k \geq 0} \mathcal{M}_{k} = \mathbb{C}[E_4(\tau), E_6(\tau)],
\end{eqnarray*}
where $E_{2k}(\tau) = 1 - \frac{2k}{B_k} \sum_{n \geq 1} \sigma_{2k-1}(n)q^n$ is the usual normalized
Eisenstein series of weight $2k \geq 4$.
The modular derivative in weight $k$ is the operator 
\begin{eqnarray*}
D_k = q\frac{d}{dq} + kE_2
\end{eqnarray*}
with $E_2(\tau) = -\frac{1}{12} + 2\sum_{n \geq 1} \sigma_1(n)q^n$. $D_k: \mathcal{M}_k \rightarrow \mathcal{M}_{k+2}$ defines a degree $2$ derivation $D: \mathcal{M} \rightarrow \mathcal{M}$ whose restriction to $\mathcal{M}_k$ is $D_k$.  We often write $D$ in place of $D_k$. $\mathcal{M}$ and $D$ generate the ring $\mathcal{R}$ of \emph{skew polynomials} whose elements are (noncommutative) \ polynomials $\sum_i m_id^i, \ m_i \in \mathcal{M},$ satisfying the relation
$dm - md = D(m)$.

\bigskip
Let the assumptions and notation be as in Section $1$.
The ring of holomorphic vector-valued modular forms with respect to $\rho$ is a $\mathbb{Z}$-graded
linear space
\begin{eqnarray*}
\mathcal{H}(\rho) = \bigoplus_{k=0}^{\infty} \mathcal{H}(k_0+2k, \rho)
\end{eqnarray*}
where $k_0$ is the least weight for which a nonzero form exists. Then
\begin{eqnarray*}
k_0 = 6(m_1+m_2)-1,
\end{eqnarray*}
and $\mathcal{H}(k_0, \rho) = \mathbb{C}F_0$
is $1$-dimensional. Let $f_1(\tau), f_2(\tau)$ be the  component functions
 of $F_0(\tau)$ with Fourier coefficients $a_{n, i}$. The leading coefficients $a_{0, i}$ are nonzero and we may, and shall, assume that $a_{0, i}=1, \ i=1, 2$. There is a natural componentwise action of $\mathcal{R}$ on elements in $\mathcal{H}(\rho)$ which turns the latter space into a $\mathbb{Z}$-graded left $\mathcal{R}$-module. Indeed, 
$\mathcal{H}(\rho) = \mathcal{R}F_0$ is a \emph{cyclic} $\mathcal{R}$-module with generator $F_0$, Moreover
$\mathcal{H}(\rho)$ is a \emph{free} $\mathcal{M}$-module with free generators $F_0, DF_0$. 

\medskip
$f_1(\tau)$ and $f_2(\tau)$ form a fundamental system of solutions of a modular linear differential equation (MLDE) of weight $k_0$ and order $2$, namely
\begin{eqnarray}\label{DE}
D_{k_0}^2f + \kappa_1E_4f = 0,
\end{eqnarray}
where
\begin{eqnarray*}
\kappa_1 = (1 - 36(m_1 - m_2)^2)/144
\end{eqnarray*}
and $D^2_{k_0} = D_{k_0+2}\circ D_{k_0}$.
The MLDE (\ref{DE})  can be rewritten as
\begin{eqnarray*}
q^2\frac{d^2f}{dq^2} + [1+2(k_0+1)E_2]q\frac{df}{dq} + [k_0(k_0+1)E_2^2 + (\kappa_1 + k_0/144)E_4]f = 0,
\end{eqnarray*}
making it clear that $q=0$ is a regular singular point.

\medskip
  There is a recursive formula for the Fourier coefficients
$a_{n, i}$ (e.g., \cite{H}, pp. 157).  We review the details as we will need them
later. Set
\begin{eqnarray*}
 \sum_{n = 0}^{\infty} u_n q^n&=& 1+2(k_0+1)E_2(q) = (1- m_1-m_2)+ O(q), \\
 \sum_{n=0}^{\infty} v_n q^n &=& k_0(k_0+1)E_2(q)^2 + (\kappa_1 + k_0/144)E_4(q)
 = m_1m_2+O(q), \\
 I_0(s) &=& s^2 - (m_1+m_2)s +m_1m_2,\\
 I_j(s) &=& u_js+v_j, \ j \geq 1.
\end{eqnarray*}
Explicitly, 
\begin{eqnarray}\label{ingres}
 u_n &=& 24(m_1+m_2)\sigma_1(n), \ n \geq 1 \notag \\
 v_n &=&k_0(m_1+m_2)(24\sum_{ r =1}^{n-1} \sigma_1(r)\sigma_1(n-r)  - 2\sigma_1(n) ) \notag \\
 && \hspace{0cm} +10(m_1+m_2-6(m_1-m_2)^2)\sigma_3(n), \ n \geq 1  \\
 I_0(n+m_1) &=&  n(n+m_1 - m_2), \notag \\
I_j(n+m_1-j) &=& (n+m_1-j)u_j + v_j, j \geq 1. \notag
\end{eqnarray}
 The recursive formula is then given (setting $a_n = a_{n, 1}$)  by $a_0 = 1$ and
\begin{eqnarray}\label{recursionforan}
a_n = - \sum_{j = 1}^n a_{n-j} \frac{I_j(m_1+n-j)}{I_0(m_1+n)}, \ \ n \geq 1.
\end{eqnarray}
There is an analogous formula for the coefficients of $f_2(\tau)$.

 \section{The projective level $M$}
We retain the notation of the previous Subsection.
 \begin{lemma}\label{lemma3.1} Suppose that $\rho$ has finite projective level $M$.
 Then $m_1, m_2$ and all Fourier coefficients $a_{n, i}$ lie in $\mathbb{Q}$.
 \end{lemma}
{\bf Proof} The assumption of the Lemma is that $\rho(T)$ has finite order $M$.
Then (\ref{rhoTdef}) shows that $m_1, m_2 \in \mathbb{Q}$. It then follows inductively
using  (\ref{ingres}) and (\ref{recursionforan}) that each Fourier coefficient is also rational. $\hfill \Box$

\bigskip
From now on we will always assume that the projective level is finite. Introduce integers
  $a, b, c, d, N$ as follows:
 \begin{eqnarray*}
&&m_1 = a/N, m_2 = b/N, \ (a, b, N) = 1, \ N > a > b \geq 0, \\
&& c = (a-b, N), a-b = cd, N = cM.
\end{eqnarray*}
$N$ is  the order of the matrix $\rho(T)$ (considered as an element of the group
$\rho(\Gamma)$), which we also call the level of $\rho$. Let
\begin{eqnarray*}
\Delta(N) = \langle \gamma T^N \gamma^{-1} \ | \ \gamma \in \Gamma \rangle
\end{eqnarray*}
be the \emph{normal closure} of $T^N$ in $\Gamma$. Thus $\Delta(N)\subseteq$ ker$\rho$
and ker$\rho$ has \emph{level $N$} in the sense of \cite{W}. 
$M$ is the projective level of $\rho$.
The following result  includes  (\ref{ab1}). 
 \begin{propn}\label{thm2.4} The following are equivalent:
\begin{eqnarray*}
&&(a)\ \rho(\Gamma) \  \mbox{is finite},\\
&&(b)\ \overline{\rho}(\Gamma) \ \mbox{is finite}, \\
&&(c) \ \rho \ \mbox{has projective level} \ M \leq 5, \\
&&(d) \ \mbox{$\rho$ is modular of level $N$}.
\end{eqnarray*}
 \end{propn}
{\bf Proof} 
 Because the index $|\Gamma: \Gamma'|$ of the commutator subgroup $\Gamma'$ in $\Gamma$ is finite, $\rho(\Gamma)$ is finite if, and only if, $\rho(\Gamma)'$ is finite. Now an old theorem of Schur says that $G'$ is a finite subgroup of a group $G$ if, and only if, the index $|G: Z(G)|$ of the center is finite. 
 Applying this to $\rho(\Gamma)$, we conclude that $\rho(\Gamma)$ is finite if, and only if,
 $|\rho(\Gamma): Z(\rho(\Gamma))|$ is finite. Since $\rho(\Gamma)$ is an irreducible group of linear transformations, we have $Z(\rho(\Gamma)) = \rho(\Gamma) \cap Z$, where we have set $Z = Z(GL(2, \mathbb{C}))$. Since
 $\overline{\rho}(\Gamma) = \rho(\Gamma)Z/Z \cong \rho(\Gamma)/Z(\rho(\Gamma))$, the equivalence of (a) and (b) in Theorem \ref{thm2.4} follows. The implication (a) $\Rightarrow$ (d) is proved in \cite{M2}, Theorem 3.5, while the converse is well-known.

 \medskip
It is also well-known (eg. \cite{W})  that if $M \leq 5$ then $\Delta(M) = \Gamma(M)$ is the principal congruence subgroup of $\Gamma$ of level $M$. Since 
$\Delta(M) \subseteq$ ker$\bar{\rho}$ it follows that if $M \leq 5$ then
ker$\overline{\rho}$ contains $\Gamma(M)$ and in particular has finite index in $\Gamma$.  This establishes the implication (c) $\Rightarrow$ (b). 
 
 \medskip
  Finally, we establish the implication (a) $\Rightarrow$ (c). Indeed, if (a) holds we may, and now shall, assume that $\rho$ is a \emph{unitary} representation. In this case,
 $\overline{\rho}(\Gamma)$ is isomorphic to a subgroup of 
 $PSU(2, \mathbb{C}) \cong SO(3, \mathbb{R})$. From Klein's classification of the finite subgroups
 of $SO(3, \mathbb{R})$, it follows that $\overline{\rho}(\Gamma)$ is one of the following:
 cyclic, dihedral, $A_4, S_4$, or $A_5$. The first case is excluded since $\rho$ is irreducible. In the last three cases, the elements of $\overline{\rho}(\Gamma)$ have order at most $5$, and in particular
 $M \leq 5$. In the dihedral case, $\overline{\rho}(T)$ generates the commutator quotient of
 $\overline{\rho}(\Gamma)$ since $T$ generates the commutator quotient of $\Gamma$, so that
$M\leq 4$ in this case. We have thus established that (c) holds in all cases, and the proof of  Theorem \ref{thm2.4} is complete. $\hfill  \Box$

\bigskip
We record some numerical restrictions on the integers $a, b, c, d, M, N$.

  \begin{lemma}\label{lemma3.21} The following hold:
 \begin{eqnarray*}
 &&(a) \ N|6(a+b) \\
 &&(b) \ (c, M)|2 \\
 &&(c) \ 4|M \Rightarrow 2||c \\
 &&(d) \ c|(N, 12) 
 \end{eqnarray*}
\end{lemma}
 {\bf Proof}: Because $\rho$ is irreducible then $\rho(S^2) = \pm I_2$, in particular $\det\rho(S^2) =1$.
It follows from the relations in $\Gamma$ that $\det\rho(T^6)=1$, that is 
$6(m_1+m_2) \in \mathbb{Z}$. Part (a) follows from this.

From (a) we obtain $c|6(a \pm b) \Rightarrow c | (12a, 12b, N) \Rightarrow c|(N, 12)$, establishing
(d).  If $3|(c, M)$ then  $3|a\pm b \Rightarrow
 3|(a, b, N) = 1$, contradiction. This shows that $(c, M)|4$. If $(c, M) = 4$ then $16|N|6(a+b)
 \Rightarrow 4|a \pm b \Rightarrow 2|(a, b, N) = 1$, a contradiction which proves (b). Finally, $4|M \Rightarrow 2|a+b
 \Rightarrow 2|a-b = cd.$ As $(d, M)=1$ then $2|c.$
  This proves (c) and completes the proof of the Lemma.  $\hfill \Box$

  \begin{lemma}\label{lemmaM6} $M \not= 6$.
   \end{lemma}
   {\bf Proof} Assume that $M=6$. Because $\rho$ is irreducible, it follows from
   \cite{M2}, Theorem 3.1 that $e^{2 \pi i (a-b)/N} = e^{2 \pi i d/6}$ is \emph{not} a primitive sixth root of unity. But then  $(d, M) = (d, 6) \not= 1$, a contradiction because $(d, N)=1$.
   $\hfill \Box$

   \begin{lemma}\label{N6812} If $N = 6, 8$ or $12$ then $M \leq 5$ and $\rho$ is modular.
   \end{lemma}
  {\bf Proof}    Suppose that $N = 6$. We are done unless perhaps $M = 6$, and this is not possible by Lemma \ref{lemmaM6}. 
     Now assume that $N = 8$ or $12$.  Since
     $N|6(a+b)$ then $2|a+b$, whence $2|a-b$ and $2|c$. If $c=2$ then either $M=4$ as required or else
     $M=6$, against Lemma \ref{lemmaM6}. On the other hand, if $c \geq 3$ then $M \leq N/3 \leq 4$. This completes the proof of the Lemma. $\hfill \Box$

 \section{Proof of Theorem 1}
 For a  prime number $p$ define $\nu_p: \mathbb{Q}^* \rightarrow \mathbb{Z}$ as follows: 
for nonzero integers $m, n,  \nu_p(n) = a$ if $n = p^al$ with $(l, p)=1$, and
$\nu_p(m/n) = \nu_p(m)-\nu_p(n)$. Set
\begin{eqnarray*}
M = (M, 60)Q.
\end{eqnarray*}

 \begin{propn}\label{thm1.4} Let $p$ be a prime such that $p|Q$. Then the $n$th Fourier coefficient
 $a_n$ of $f_1(\tau)$ satisfies
\begin{eqnarray}\label{nuform}
\nu_p(a_n) =  - \nu_p(Q^nn!).
\end{eqnarray}
In particular, $\nu_p(a_n)$ is \emph{strictly decreasing} for 
$n \geq 0$.
 \end{propn}
 {\bf Proof}:  
 We prove the Theorem by induction on $n$. The case $n=0$ holds because $a_0=1$. In the notation of Section $3$, the recursion  (\ref{recursionforan}) reads as follows:
\begin{eqnarray}\label{yetmorean}
a_n &=&  \frac{-1}{nc(Mn+d)} \sum_{j = 1}^n a_{n-j}\left\{(a +N(n-j))u_j+Nv_j)\right\} \notag \\
&=& \frac{-1}{nc(Mn+d)} \sum_{j = 1}^n a_{n-j}s_j(n), 
\end{eqnarray}
where we have set 
\begin{eqnarray*}
s_j(n) =  (a+N(n-j))u_j+Nv_j.
\end{eqnarray*}
Then
\begin{eqnarray}\label{sdef}
  &&\hspace{4cm} s_j(n) =  \notag \\
 &=& \left\{\frac{2(a+b)\sigma_1(j)( (12n-12j+1)N+6(a-b))}{N} \right\} + 
  \frac{10\sigma_3(j) ((a+b)N-6(a-b)^2)}{N} \notag \\
 &&\hspace{0.3cm} + \frac{24(a+b)(6(a+b)-N)}{N}  \sum_{ r =1}^{j-1}  \sigma_1(r)\sigma_1(j-r) \notag \\
 &\equiv&  -  \frac{60\sigma_3(j) cd^2}{M} \ (\mbox{mod} \ 2\mathbb{Z}),
\end{eqnarray}
where we used Lemma \ref{lemma3.21}(a) for the last congruence.
Using Lemma \ref{lemma3.21}(b), we deduce from (\ref{sdef}) that  
\begin{eqnarray}\label{nuineq}
\nu_p\left( \frac{s_j(n)}{c}\right) \geq  - \nu_p(Q), 
\end{eqnarray}
with equality if $j=1$.

\bigskip
Let 
\begin{eqnarray}\label{bnjdef}
b_{n, j} = \frac{a_{n-j}}{n}\frac{1}{Mn+d}\frac{s_j(n)}{c} 
\end{eqnarray}
denote the $j$th. term on the right hand side of (\ref{yetmorean}). 
 Using induction,  (\ref{nuineq}) and $(M, d)=1$,  we  have for $1\leq j \leq n$ that
  \begin{eqnarray*}
 \nu_p(b_{n,j}) &\geq&  - \nu_p(Q^{n-j}(n-j)!) - \nu_p(n) - \nu_p(Q) \geq - \nu_p(Q^nn!),
 \end{eqnarray*}
and the two inequalities are both equalities if, and only if, $j=1$.
By the non-Archimedian property of $p$-adic valuations, it follows that 
\begin{eqnarray*}
\nu_p(a_n) = \nu_p \left(\sum_{j=1}^n b_{n, j} \right) = \nu_p(b_{n, 1}) =  -\nu_p(Q^nn!).
\end{eqnarray*}
This completes the proof of the Proposition. $\hfill \Box$

\bigskip
 Let $\mathcal{M}_{\mathbb{Q}}$ be the space of (classical) holomorphic modular forms with rational Fourier coefficients. Similarly, 
$\mathcal{H}(\rho)_{\mathbb{Q}}$ and $\mathcal{H}(k, \rho)_{\mathbb{Q}}$ are the corresponding spaces of vector-valued modular forms whose components have rational Fourier coefficients.
Then $\mathcal{H}(\rho)_{\mathbb{Q}}$ is a left 
$\mathcal{M}_{\mathbb{Q}}$-module, and $F_0 \in \mathcal{H}(k_0, \rho)_{\mathbb{Q}}$ by Lemma \ref{lemma3.1}.

\begin{lemma}\label{lemma5.0} The following hold:
\begin{eqnarray*}
&&(a) \ \mathcal{H}(\rho)_{\mathbb{Q}} =  \{\alpha F_0 + \beta DF_0 \ | \ \alpha, \beta \in 
\mathcal{M}_{\mathbb{Q}} \} \\
&&(b) \ \  \ \mathcal{H}(\rho) = \mathbb{C} \otimes_{\mathbb{Q}} \mathcal{H}(\rho)_{\mathbb{Q}} 
\end{eqnarray*}
\end{lemma}
{\bf Proof}: Since $F_0 \in \mathcal{H}(\rho)_{\mathbb{Q}}$ then also $DF_0 \in  \mathcal{H}(\rho)_{\mathbb{Q}}$, and we conclude that  $\{\alpha F_0 + \beta DF_0 \ | \ \alpha, \beta \in \mathcal{M}_{\mathbb{Q}} \} \subseteq \mathcal{H}(\rho)_{\mathbb{Q}}$. Because $\mathcal{M} = \mathbb{C} \otimes_{\mathbb{Q}} \mathcal{M}_{\mathbb{Q}}$ and $\mathcal{H}(\rho) = \mathcal{M}F_0+\mathcal{M}DF_0$, part (b) follows. 

\medskip
Now let $G \in \mathcal{H}(\rho)_{\mathbb{Q}}$. Because $F_0, DF_0$ are free generators
of $\mathcal{H}(\rho)$ considered as $\mathcal{M}$-module, there are unique forms 
$\alpha, \beta \in \mathcal{M}$ with the property that
   $\alpha F_0 + \beta DF_0 = G$. To complete the proof of part (a) we have to show that $\alpha, \beta \in \mathcal{M}_{\mathbb{Q}}$. 
   
   \medskip
   Let
$g_i = \sum_{n \geq 0} b_{n, i}q^{n+m_i}, i = 1, 2$ be the components of $DF_0$. Thus, each $a_{n, i}, b_{n, i} \in \mathbb{Q}$.
Let $\alpha = \sum_{n \geq 0} c_nq^n, \beta = \sum_{n \geq 0} d_nq^n$. 
We prove by induction on $n$ that each $c_n, d_n \in \mathbb{Q}$. 

\medskip
The leading column of $F_0$
is $(q^{m_1}, q^{m_2})^t$, while that for $DF_0$ is \linebreak 
$((m_1 - k_0)/12q^{m_1}, (m_2 - k_0/12)q^{m_2})^t$. Thus the $n$th coefficient 
 column for $G$ is equal to 
 \begin{eqnarray}\label{colcoeffs}
c_n \left(\begin{array}{c}1 \\1\end{array}\right) + \sum_{j=1}^n c_{n-j} \left(\begin{array}{c}a_{j, 1} \\a_{j, 2}\end{array}\right) + d_n\left(\begin{array}{c} m_1-k_0/12 \\ m_2 - k_0/12\end{array}\right)
+  \sum_{j=1}^n d_{n-j} \left(\begin{array}{c}b_{j, 1} \\b_{j, 2}\end{array}\right) 
\end{eqnarray}
and by hypothesis this has rational entries. If $n=0$ this says that
 \begin{eqnarray*}
c_0 \left(\begin{array}{c}1 \\1\end{array}\right)  + d_0\left(\begin{array}{c} m_1-k_0/12 \\ m_2 - k_0/12\end{array}\right)
\end{eqnarray*}
is rational. We can write this in the form
\begin{eqnarray*}
\left(\begin{array}{cc} 1 & m_1 - k_0/12 \\  1 & m_2 - k_0/12\end{array}\right) \left(\begin{array}{c} c_0 \\ d_0\end{array}\right) \in \mathbb{Q}^2.
\end{eqnarray*}
Since $m_1 \not= m_2$ the matrix in the last display is invertible, so $c_0, d_0 \in \mathbb{Q}$. This begins the induction.
The inductive step follows by the same argument, using (\ref{colcoeffs}). This completes the proof of the Lemma.  $\hfill \Box$

  \bigskip
 For a prime $p$, let $\mathcal{B}_p \subseteq \mathcal{H}(\rho)_{\mathbb{Q}}$ be the span of those vector-valued modular forms with the property
 that  there is an integer $B$ such that $BF(\tau)$ has $p$-integral Fourier coefficients.
  It is well-known that the Fourier coefficients of every element in $\mathcal{M}_{\mathbb{Q}}$ have this property.
Consequently, $\mathcal{B}_p$ is an $\mathcal{M}_{\mathbb{Q}}$-submodule of $\mathcal{H}(\rho)_{\mathbb{Q}}$.

     \begin{propn}\label{theorem4.3}   Suppose that $p|Q$. Then $\mathcal{B}_p = \{0\}$. 
\end{propn}
{\bf Proof}:  We assume that $ \mathcal{B}_p \not= \{0\}$ and derive a contradiction. We first prove
\begin{eqnarray}\label{spec}
  \mathcal{M}_{\mathbb{Q}}F_0  \cap \mathcal{B}_p = \{0\}.
\end{eqnarray}
 
 Because $\nu_p(Q) \geq 1$, we know  that (\ref{nuform}) holds. 
Suppose that  $0 \not= \alpha \in \mathcal{M}_{\mathbb{Q}}$ with $\alpha = \sum_{n\geq 0} b_nq^n$. There is an integer $N$ such that $\nu_p(b_n) \geq N$ for all $n$. Let $c_n$ be the $n$th. Fourier coefficient of the first component of $\alpha F_0$, so that
 \begin{eqnarray}\label{cnsum}
c_n = \sum_{j=0}^n b_{n-j}a_j.
\end{eqnarray}
Assume to begin with that $b_0 \not= 0$, and choose any integer $s >\nu_p(b_0)-N$.
 Using (\ref{nuform})
we find that if $j < p^s$ then either $b_{p^s-j}a_j = 0$ or else
\begin{eqnarray*}
&&\hspace{3.5cm} \nu_p(b_0 a_{p^s}) -\nu_p(b_{p^s-j}a_j)  \\
&=& \nu_p(b_0)  - p^s\nu_p(Q) - \nu_p(p^s!)
-\nu_p(b_{p^s-j})  + j\nu_p(Q) + \nu_p(j!) \\
&\leq& \nu_p(b_0) -\nu_p(Q)(p^s-j)  - s - N \\
&<& \nu_p(b_0)  - s - N < 0.
\end{eqnarray*}
Now (\ref{cnsum}) and the non-Archimedian property of $\nu_p$ shows that
\begin{eqnarray*}
\nu_p(c_{p^s}) = \nu_p(b_0a_{p^s}) &=&  \nu_p(b_0) - p^s\nu_p(Q) - \nu_p(p^s!) \\
&\leq&  \nu_p(b_0) - p^s- s.
\end{eqnarray*}
Therefore $\alpha F_0$ does \emph{not} have bounded $p$-power.
If $b_0 = 0$ we obtain a similar result by using the first nonvanishing coefficient of
$\alpha$ in place of $b_0$ in the previous argument. This completes the proof of (\ref{spec}).

\medskip
 Now let $0 \not= F \in \mathcal{B}_p$ have weight $k$. By Lemma \ref{lemma5.0}(a) we have 
 \begin{eqnarray}\label{Feq}
 F =\alpha F_0 + \beta DF_0
 \end{eqnarray}
 with $\alpha, \beta \in \mathcal{M}_{\mathbb{Q}}$.  Note that $\beta \not= 0$ by (\ref{spec}).

 \medskip
  Apply $D$ to (\ref{Feq})  to get
  \begin{eqnarray*}
(D\alpha)F_0 + \alpha DF_0 + (D\beta)DF_0 + \beta D^2F_0 = DF,
\end{eqnarray*}
 and then use (\ref{DE})  to obtain
\begin{eqnarray}\label{DE2}
(D\alpha) F_0 + \alpha  DF_0 + (D\beta)DF_0 -\kappa_1 \beta E_4 F_0  = DF.
\end{eqnarray}
Now use (\ref{Feq}) and (\ref{DE2}) to eliminate the $DF_0$ terms. We obtain
    \begin{eqnarray}\label{fe}
   ((D\alpha)\beta  -\kappa_1 \beta^2 E_4 -\alpha^2 -\alpha D\beta)F_0 = \beta DF -(\alpha + D\beta)F.
   \end{eqnarray}
The coefficient of $F_0$ on the left-hand-side of (\ref{fe}) lies in $\mathcal{M}_{\mathbb{Q}}$, while
the right-hand-side lies in $\mathcal{B}_p$. Thanks to (\ref{spec}), the only way this can happen is
for both sides to be identically zero. Then the components of $F$ satisfy the first order differential equation
 \begin{eqnarray}\label{fe}
 \beta DF -(\alpha + D\beta)F = 0, 
   \end{eqnarray}
and since they are linearly independent this is not possible. This completes the proof of Proposition \ref{theorem4.3}, and with it that of Theorem 1 as well. $\hfill \Box$

\medskip It follows easily from Lemma \ref{lemma3.21} that if $M|60$ then $N|120$. So the version of the main Theorem stated in the abstract is indeed a special case of Theorem 1.

\section{Concluding remarks}
Suppose that $\rho$ is \emph{not} modular, and that $F(\tau) \in \mathcal{H}(k, \rho)_{\mathbb{Q}}$
is a nonzero vector-valued modular form with bounded denominators. 
From Theorem 1 and Lemmas \ref{lemmaM6} and \ref{N6812} it follows that
$M|60, N|120$. Indeed, there are only $13$ possible pairs $(M, N)$, namely
\begin{eqnarray}\label{pairs}
&&(10, 10), (10, 20), (10, 30), (10, 60),  (12, 24), (15, 15), (15, 30), \notag\\
  && (15, 60), (20, 40),
 (20, 120), (30, 30), (30, 60), (60, 120).
  \end{eqnarray}
It is straightforward, though unenlightening, to enumerate all of  the equivalence classes of $\rho$ corresponding to these values using Theorem 3.1 of \cite{M2}. It turns out that there are 282 such classes.

\medskip
The recursive formula (\ref{recursionforan}) is convenient for  machine calculation. Inspection of  the denominators of the first thousand or so coefficients of  the  examples corresponding to the pairs $(M, N)$ in (\ref{pairs}) shows that they steadily increase, with more and more primes occurring in the denominators. There seems to be  no question that the Conjecture is true! On the other hand,  it is unclear whether there are primes $p$ (analogous to those dividing $Q$ in Proposition
\ref{thm1.4}) which divide the denominators of the Fourier coefficients of $f_1$ to a power which is unbounded for $n \rightarrow \infty$.  It is also unclear whether one should expect that all but a finite number of primes occur in the denominators of these Fourier coefficients.
This property holds for certain classes of generalized modular forms (\cite{KoM1}, \cite{KoM2}) and is something that one expects to be a rather general feature of vector-valued modular forms. The arithmetic nature of the denominators of the Fourier coefficients of the components of $F_0$ - whether $M$ divides $60$ or not - appears to be quite interesting.

 \bigskip
 \noindent
 Author's address: Department of Mathematics,  University of California,  Santa Cruz, 
CA 95064, U.S.A. 
 E-mail: gem@cats.ucsc.edu


\begin{thebibliography}{BPZ}
   
   \bibitem[AS]{AS} Atkin, O., and Swinnerton-Dyer, P., Modular forms on noncongruence subgroups,
   \textit{Proc. Symp. Pure Math. Vol. XIX}, Amer. Math. Soc., Providence, R.I. (1971), 1-25.
 
 \bibitem[H]{H} Hille, E., \textit{Ordinary Differential Equations in the
 Complex Domain}, Dover Publications, New York, 1976. 
 
 \bibitem[KM]{KM} Knopp, M., and Mason, G., On vector-valued modular forms and their Fourier coefficients, Acta Arithmetica \textbf{110} No. 2 (20033), 117-125.
 
 \bibitem[KoM1]{KoM1} Kohnen, W. and Mason, G., On Generalized Modular Forms and their Applications,
 Nagoya J. Math. \textbf{192} (2008),  119-136.

   \bibitem[KoM2]{KoM2} Kohnen, W. and Mason, G., On the canonical decomposition of a generalized modular form, submitted, (Arxiv: 1003.2407.)
 
 \bibitem[KL1]{KL1} Kurth, C. and Long, L., On modular forms for some noncongruence subgroups
 of $SL_2(\mathbb{Z})$, J. Numb. Th. \textbf{128} No. 7 (2008), 1989-2009.
 
 \bibitem[KL2]{KL2} Kurth, C. and Long, L., On modular forms for some noncongruence subgroups
 of $SL_2(\mathbb{Z})$ II, Bull. Lond. Math. Soc. \textbf{41} No. 4 (2009), 589-598.

  \bibitem[MM]{MM} Marks, C., and Mason, G., Structure of the module of vector-valued modular forms,
  J. Lond. Math. Soc. \textbf{82} Part $1$ (2010), 32-48.
  
\bibitem[M1]{M1} Mason, G., Vector-valued modular forms and linear differential equations, 
Int. J. Numb. Th. \textbf{3} No. 3 (2007), 1-14.

\bibitem[M2]{M2} Mason, G., 2-dimensional vector-valued modular forms,
 Ramanujan Journal \textbf{17} (2008), 405-427.

 \bibitem[W]{W} Wohlfahrt, K., An extension of F. Klein's level concept, Ill. J. Math. \textbf{13}
 (1964), 529-535.
 
  \end{thebibliography}
\end{document}